
\baselineskip=14pt
\parskip=10pt

\magnification=\magstephalf

\def\B{{\cal B}}

\def\1{{\overline{1}}}
\def\2{{\overline{2}}}
\parindent=0pt
\overfullrule=0in

\def\frac#1#2{{#1 \over #2}}
\centerline
{\bf 
Explicit Expressions for the First 20 Moments of the Area Under Dyck and Motzkin Paths
}
\bigskip
\centerline
{\it AJ BU, Shalosh B. EKHAD, and Doron ZEILBERGER}
\bigskip

{\bf Abstract}: 
We show the utility of AJ Bu's recent article for computing explicit expressions for the {\it generating functions} of sums of powers of areas 
under Dyck and Motzkin paths, by deducing from them {\bf explicit} expressions for the {\it actual sequences}. 
This leads to explicit expressions for the scaled moments, and to the elementary verification, up to any specific moment,
that they tend to those of the area under Brownian excursion.

{\bf Maple package and Output Files}

This article is accompanied by a Maple package, {\tt qEWplus.txt}, and two output files.

They are obtainable from the front of this article

{\tt https://sites.math.rutgers.edu/\~{}zeilberg/mamarim/mamarimhtml/qew.html } \quad .

{\bf The Power of Symbolic Computation}

In [B], the first author harnessed the power of symbolic computation to generate {\it generating functions} for the sum of $r^{th}$ powers of areas under Dyck and Motzkin paths, by
interfacing algebra and calculus. Once you have these sequences themselves, you can immediately derive the actual moments (by dividing by the
total number of paths, i.e. the Catalan and Motzkin numbers respectively), that in turn would give you the {\it moments about the mean} and the {\it scaled moments}.
Alas, getting {\it generating functions} for these quantities is awkward. It would be much more helpful to have {\bf explicit expressions} for the {\it actual} sequences,
not just their generating functions.

{\bf The power of guessing}

Recall that a {\it Dyck path} (respectively {\it Motzkin path}) is a sequence of {\it atomic steps}  of the form
$(x,y) \rightarrow (x+1,y-1),(x,y) \rightarrow (x+1,y+1)$ (respectively of the form $(x,y) \rightarrow (x+1,y-1),(x,y) \rightarrow (x+1,y+1), (x,y) \rightarrow (x+1,y)$)
that

$\bullet$ Start at the origin $(0,0)$.

$\bullet$ End somewhere on the  $x$ axis, i.e. $y=0$ .

$\bullet$ Never go strictly below the  $x$ axis, i.e. always stay in $y \geq 0$.

Equivalently sequences in $\{1,-1\}$ (respectively $\{-1,0,1\}$) whose total sum is $0$ and whose partial sums are non-negative. See [B] for details.

It turns out (the proofs are left to the reader) that the following facts hold.

{\bf Fact 1}: Let $C_r(n)$ be the sum of the $r^{th}$ powers of the areas under Dyck paths from $(0,0)$ to $(2n,0)$:
$$
C_r(n):= \sum_{P \in Dyck(n)} \, Area(P)^r \quad .
$$
Also let $C_n$ be the Catalan numbers $(2n)!/(n!(n+1)!)$ (same as $C_0(n)$).

There exist {\it polynomials} in $n$, $P_1(n)$ and  $P_2(n)$ such that
$$
C_r(n)= P_1(n) \, C_n \,+ \,P_2(n) \,4^n \quad .
$$

{\bf Fact 2}: Let $M_r(n)$ be the sum of the $r^{th}$ powers of the areas under Motzkin paths from $(0,0)$ to $(n,0)$:
$$
M_r(n):= \sum_{P \in Motzkin(n)} \, Area(P)^r \quad .
$$
Also let $M_n$ be the Motzkin numbers $\sum_{r=0}^{n}  \frac{n! (2r)!}{r!^2(r+1)!(n-r)!}$ (same as $M_0(n)$).

There exist {\it polynomials} $Q_1(n)$, $Q_2(n)$ ,  $Q_3(n)$ , $Q_4(n)$  such that
$$
M_r(n)= Q_1(n) \, M_n \,+ \, Q_2(n) \, M_{n+1} \,+ \, Q_3(n) \, 3^n \, + \, Q_4(n) \, (-1)^n  \quad .
$$

Since we know the {\it format} of these closed-form expressions for $C_r(n)$ and $M_r(n)$, we can crank out sufficiently many terms, use {\it undetermined coefficients},
and ask Maple to {\tt solve} the resulting system of linear equations. Note that {\bf once conjectured}, the proof is {\it purely routine}, even if you don't
believe that the facts always hold. Maple can {\it all by itself} find the generating functions for the conjectured expressions (they all are derivable from
those for the Catalan and Motzkin numbers), and check that they coincide with those found in [B].

Once we have explicit expressions for $C_r(n)$, we get an explicit expression for the $r^{th}$-moment $m_r(n)=\frac{C_r(n)}{C_n}$, and from them the {\it moments about the mean}, and in turn,
the scaled moments. 
Similarly for the Motzkin numbers (and any of their ilk). 
Then we take the limit as $n \rightarrow \infty$.

To our great delight, they coincided, up to the $20$-th moment with those of the area under Brownian Excursion $\B_{ex} ([T]$,  and  see [J] for a great exposition). 
Of course one can go much further.

According to [J] (Equations $(6)$ and $(7)$ there), 
$$
{\bf E} \B^{k}_{ex} \, = \, \frac{4 \sqrt{\pi} 2^{-k/2}k!}{\Gamma((3k-1)/2)} \, K_k \quad ,  \quad k \geq 0 \quad,
$$
where $K_k$ is the sequence of rational numbers defined by the recurrence, $K_0=-\frac{1}{2}$ and
$$
K_k = \frac{3k-4}{4} K_{k-1} + \sum_{j=1}^{k-1} K_j \, K_{k-j} \quad, \quad k \geq 0 \quad.
$$

Once you centralize and scale them, you get {\bf exactly} what we got for the corresponding limits of the first (centralized-scaled) $20$ moments for the Dyck paths and Motzkin paths, by fully {\it elementary} and {\it finitistic} means.

But we can do much more than {\it just} find the {\it limit}. Since we have explicit expressions, we can take the {\it asymptotics} to {\it any} desired order.

{\bf Explicit expressions for the sum of the powers of areas under Dyck paths from $(0,0)$ to $(2n.0)$ for all powers from $1$ to $20$.}

$$
C_1(n)=\frac{\left(2 n \right)! \left(-2 n -1\right)}{n ! \left(1+n \right)!}+4^{n} \quad,
$$

$$
C_2(n)=\frac{\left(2 n \right)! \left(\frac{10}{3} n^{3}+11 n^{2}+\frac{26}{3} n +2\right)}{n ! \left(1+n \right)!}+4^{n} \left(-4 n -2\right) \quad,
$$

$$
C_3(n)=\frac{\left(2 n \right)! \left(-20 n^{4}-60 n^{3}-61 n^{2}-26 n -4\right)}{n ! \left(1+n \right)!}+4^{n} \left(\frac{15}{4} n^{3}+\frac{75}{4} n^{2}+\frac{33}{2} n +4\right) \quad,
$$

$$
C_4(n)=
\frac{\left(2 n \right)! \left(\frac{884}{63} n^{6}+\frac{4568}{35} n^{5}+\frac{108043}{315} n^{4}+\frac{14016}{35} n^{3}+\frac{74332}{315} n^{2}+\frac{2416}{35} n +8\right)}{n ! \left(1+n \right)!}
$$
$$
+4^{n} \left(-30 n^{4}-101 n^{3}-111 n^{2}-50 n -8\right) \quad ,
$$

$$
C_5(n)=
\frac{\left(2 n \right)! \left(-\frac{8840}{63} n^{7}-\frac{53044}{63} n^{6}-\frac{128174}{63} n^{5}-\frac{162091}{63} n^{4}-\frac{116888}{63} n^{3}-\frac{48520}{63} n^{2}-\frac{3608}{21} n -16\right)}{n ! \left(1+n \right)!}
$$
$$
+4^{n} \left(\frac{565}{32} n^{6}+\frac{6495}{32} n^{5}+\frac{18485}{32} n^{4}+\frac{22825}{32} n^{3}+\frac{7035}{16} n^{2}+\frac{535}{4} n +16\right) \quad .
$$

For the explicit expressions for $C_r(n)$ for $6 \leq r \leq 20$, see the output file:

{\tt http://sites.math.rutgers.edu/\~{}zeilberg/tokhniot/oqEWplus1.txt} \quad .

{\bf Explicit expressions for the sum of the powers of areas under Motzkin paths from $(0,0)$ to $(n,0)$ for all powers from $1$ to $20$.}

As above $M_n$ (same as $M_0(n)$) denote the Motzkin numbers.

$$
M_1(n)=\left(-n -1\right) M_{n}+\frac{3^{n+1}}{4}+\frac{\left(-1\right)^{n}}{4} \quad,
$$

$$
M_2(n)=
$$
$$
\left(\frac{25}{36} n^{3}+\frac{11}{3} n^{2}+\frac{191}{36} n +\frac{7}{3}\right) M_{n}
+ \left(-\frac{5}{36} n^{3}-\frac{2}{3} n^{2}-\frac{31}{36} n -\frac{1}{3}\right)  M_{n+1}
+\left(-\frac{3 n}{2}-\frac{3}{2}\right) 3^{n}+\left(-\frac{n}{2}-\frac{1}{2}\right) \left(-1\right)^{n} \quad,
$$

$$
M_3(n)=
$$
$$
\left(-\frac{25}{12} n^{4}-\frac{133}{12} n^{3}-\frac{251}{12} n^{2}-\frac{203}{12} n -5\right) M_{n}+ \left(\frac{5}{12} n^{4}+\frac{29}{12} n^{3}+\frac{55}{12} n^{2}+\frac{43}{12} n +1\right) M_{n+1}
$$
$$
+\left(\frac{15}{64} n^{3}+\frac{207}{64} n^{2}+\frac{1521}{256} n +\frac{1503}{512}\right) 3^{n}
+\left(\frac{15}{64} n^{3}+\frac{93}{64} n^{2}+\frac{585}{256} n +\frac{545}{512}\right) \left(-1\right)^{n} \quad,
$$

$$
M_4(n)=
\left(\frac{221}{324} n^{6}+\frac{31441}{3780} n^{5}+\frac{411613}{11340} n^{4}+\frac{72251}{945} n^{3}+\frac{34189}{405} n^{2}+\frac{59401}{1260} n +\frac{3314}{315}\right) M_{n}
$$
$$
+ \left(-\frac{221}{1134} n^{6}-\frac{1409}{630} n^{5}-\frac{27248}{2835} n^{4}-\frac{10781}{540} n^{3}-\frac{122653}{5670} n^{2}-\frac{6349}{540} n -\frac{794}{315}\right) M_{n+1}
$$
$$
+\left(-\frac{15}{16} n^{4}-\frac{63}{8} n^{3}-\frac{1197}{64} n^{2}-\frac{2241}{128} n -\frac{735}{128}\right) 3^{n}+\left(-\frac{15}{16} n^{4}-\frac{19}{4} n^{3}-\frac{573}{64} n^{2}-\frac{947}{128} n -\frac{289}{128}\right) \left(-1\right)^{n}
\quad .
$$

For the explicit expressions  of $M_r(n)$, for $5 \leq r \leq 20$, see the output file:

{\tt https://sites.math.rutgers.edu/\~{}zeilberg/tokhniot/oqEWplus2.txt} \quad .

{\bf Acknowledgment}: Many thanks to Svante Janson for educating us on Brownian excursions.

{\bf References}

[B] AJ Bu,  {\it Explicit generating functions for the sum of the areas under Dyck and Motzkin paths and for their powers}, to appear in
Discrete Mathematics Letters. \hfill\break
arxiv: {\tt https://arxiv.org/abs/2310.17026} \quad ; \hfill\break
\quad  front: {\tt https://ajbu1.github.io/Papers/MotzArea/MotzArea.html} \quad .

[J] Svante Janson, {\it Brownian excursion area, Wright's constants in graph enumeration, and other Brownian areas}. Probability Surveys {\bf 3} (2007), 80-145. \hfill\break
{\tt https://www2.math.uu.se/\~{}svante/papers/sj201.pdf}

[T] Lajos Takacs, {\it A Bernoulli excursion and its various applications}, Adv. in Appl. Probab. {\bf 23} (1991), 557-585.

\vfill\eject

\bigskip
\hrule
\bigskip
AJ Bu, Shalosh B. Ekhad and Doron Zeilberger, Department of Mathematics, Rutgers University (New Brunswick), Hill Center-Busch Campus, 110 Frelinghuysen
Rd., Piscataway, NJ 08854-8019, USA. \hfill\break
Email: {\tt    aj dot bu131 at gmail dot com} \quad,  \quad {\tt ShaloshBEkhad at gmail dot com}, \quad {\tt DoronZeil at gmail dot com}   \quad .

Written: {\bf May  5, 2024}.

{\bf Exclusively published in the Personal Journal of Shalosh B. Ekhad and Doron Zeilerger, AJ Bu's web-site, and arxiv.org} \quad .

\end